\documentclass[10pt]{amsart}
\usepackage{hyperref}
\usepackage{float}
\usepackage{citeref}
\usepackage{amscd,amssymb,amsmath,latexsym,bm}
\usepackage[mathcal,mathscr]{euscript}
\usepackage{lipsum}
\usepackage{amsfonts}
\usepackage{graphicx}
\usepackage{epstopdf}
\usepackage{algorithmic}
\usepackage{hyperref}
\usepackage{xcolor}
\usepackage{upgreek}
\usepackage{enumerate}
\usepackage{ulem}
\usepackage{gensymb}			
\usepackage[utf8]{inputenc}
\usepackage{tikz}
\usepackage{braids}

\newtheorem{code}{Computer Code}

\numberwithin{equation}{section}
\numberwithin{figure}{section}

\newcommand{\NM}{{\mathbb N}}

\newcommand{\Bb}{{\mathcal B}}

\tikzset{
  pt/.style={insert path={node[scale=2]{.}}},
  dnup/.style={insert path={ [pt] .. controls +(0,0.5) and +(0,-0.5) .. +(#1,1) [pt]}},
  dndn/.style={insert path={ [pt] .. controls +(0,0.5) and +(0,0.5) .. +(#1,0) [pt]}},
  upup/.style={insert path={ [pt] .. controls +(0,-0.5) and +(0,-0.5) .. +(#1,0) [pt]}},
}

\begin{document}

\title{Algorithm for computing Representations of the Braid Group and Temperley-Lieb algebra}

\author{Yitzchak Shmalo}

\address{Department of Mathematics 
\\Penn State, 
\\State College, PA 16801, USA \\
\href{mailto:yms5281@psu.edu}{yms5281@psu.edu}}

\thanks{I would like to thank Emil Prodan for introducing me to braid groups, useful discussions and help with this work.}

\date{\today}

\begin{abstract} The braid group appears in many scientific fields and its representations are instrumental in understanding topological quantum algorithms, topological entropy, classification of manifolds and so on. In this work, we study planer diagrams which are Kauffman's  reduction of the braid group algebra to the Temperley-Lieb algebra. We introduce an algorithm for computing all planer diagrams in a given dimension. The algorithm can also be used to multiply planer diagrams and find their matrix representation.   
\end{abstract}

\maketitle

\setcounter{tocdepth}{2}
{\scriptsize
\tableofcontents
}

\section{Introduction to Braid Group} 

\subsection{Braid Group} There are many good textbooks on the braid group and its representations. Among them, the one by Kassel and Tuarev stands out \cite{BraidBook1}.

The braid group $\Bb_N$, as a set, consists of all finite braid operations that can be applied to $N$ straight strands \cite{BraidBook1}. A braid operation can be represented pictorially by its action on $N$ straight strands, which results in $N$ weaved strands like the one shown below:
$$
\begin{tikzpicture}
\braid[number of strands=4,width=0.4cm,height=0.4cm,border height=0.1cm]  a_1 a_2^{-1} a_3 a_2^{-1};
\end{tikzpicture}
$$
In these pictures, we will take the view that time progresses vertically up. When operations are performed, it is useful to enclose the braid operations in boxes. For example, two braid operations applied sequentially in time result in a new overall braid operation, and this is better represented with boxes, as below:
$$
\begin{matrix}
\begin{tikzpicture}
\braid[number of strands=4,width=0.4cm,height=0.4cm,border height=0.1cm]   a_1  a_2^{-1};
\draw (0.2,-1) rectangle (1.8,-0);
\end{tikzpicture}
\end{matrix}\
\begin{matrix}
\cdot
\end{matrix}\
\begin{matrix}
\begin{tikzpicture}
\braid[number of strands=4,width=0.4cm,height=0.4cm,border height=0.1cm]  a_3  a_2^{-1};
\draw (0.2,-1) rectangle (1.8,-0);
\end{tikzpicture}
\end{matrix}\
\begin{matrix}
=
\end{matrix}\
\begin{matrix}
\begin{tikzpicture}
\braid[number of strands=4,width=0.4cm,height=0.4cm,border height=0.1cm,fill=] a_1 a_2^{-1} a_3 a_2^{-1};
\draw (0.2,-1.8) rectangle (1.8,-0.9);
\draw (0.2,-0.9) rectangle (1.8,0);
\end{tikzpicture}
\end{matrix}\
\begin{matrix}
=
\end{matrix}\
\begin{matrix}
\begin{tikzpicture}
\braid[number of strands=4,width=0.4cm,height=0.4cm,border height=0.1cm,fill=] a_1 a_2^{-1} a_3 a_2^{-1};
\draw (0.2,-1.8) rectangle (1.8,0);
\end{tikzpicture}
\end{matrix}
$$
In fact, the above rule defines the braid group composition. 

\vspace{0.1cm}

The braid operations are defined up to the Reidemeister moves \cite{KauffmanBook1}. For example:
\begin{equation}\label{Eq:121}
\begin{matrix} 
\begin{tikzpicture}
\braid[number of strands=3,width=0.4cm,height=0.4cm,border height=0.1cm] a_1;
\draw (0.2,-0.6) rectangle (1.4,0);
\end{tikzpicture}
\end{matrix}
\
\begin{matrix} 
\cdot
 \end{matrix}
 \
 \begin{matrix} 
\begin{tikzpicture}
\braid[number of strands=3,width=0.4cm,height=0.4cm,border height=0.1cm] a_2;
\draw (0.2,-0.6) rectangle (1.4,0);
\end{tikzpicture} 
\end{matrix}
\
\begin{matrix} 
\cdot
\end{matrix}
\
\begin{matrix} 
\begin{tikzpicture}
\braid[number of strands=3,width=0.4cm,height=0.4cm,border height=0.1cm] a_1;
\draw (0.2,-0.6) rectangle (1.4,0);
\end{tikzpicture} 
\end{matrix}
\
\begin{matrix} 
=
\end{matrix}
\
\begin{matrix} 
\begin{tikzpicture}
\braid[number of strands=3,width=0.4cm,height=0.3cm,border height=0.1cm]  a_1  a_2  a_1 ;
\draw (0.2,-1.1) rectangle (1.4,0);
\end{tikzpicture}
\end{matrix}
\end{equation}
and:
\begin{equation}\label{Eq:212}
\begin{matrix} 
\begin{tikzpicture}
\braid[number of strands=3,width=0.4cm,height=0.4cm,border height=0.1cm] a_2;
\draw (0.2,-0.6) rectangle (1.4,0);
\end{tikzpicture}
\end{matrix}
\
\begin{matrix} 
\cdot
 \end{matrix}
 \
 \begin{matrix} 
\begin{tikzpicture}
\braid[number of strands=3,width=0.4cm,height=0.4cm,border height=0.1cm] a_1;
\draw (0.2,-0.6) rectangle (1.4,0);
\end{tikzpicture} 
\end{matrix}
\
\begin{matrix} 
\cdot
\end{matrix}
\
\begin{matrix} 
\begin{tikzpicture}
\braid[number of strands=3,width=0.4cm,height=0.4cm,border height=0.1cm] a_2;
\draw (0.2,-0.6) rectangle (1.4,0);
\end{tikzpicture} 
\end{matrix}
\
\begin{matrix} 
=
\end{matrix}
\
\begin{matrix} 
\begin{tikzpicture}
\braid[number of strands=3,width=0.4cm,height=0.3cm,border height=0.1cm]  a_2  a_1  a_2 ;
\draw (0.2,-1.1) rectangle (1.4,0);
\end{tikzpicture}
\end{matrix}
\end{equation}
are equivalent from a topological point of view and they are considered as representing the same element of the braid group. Note that the ends of the strands are kept fixed when deformations of the strands are attempted. Diagrams like the ones above are called links and will generically be referred to and represented as $\boxed{\ \ L \ \ }$. Such diagram represents the entire equivalence class of that link relative to deformations allowed by topology. One can think of the explicit link being drawn as a representative and of the link plus the box as the entire equivalence class. The composition rule says that, when multiplying two elements of the braid group, one chooses representatives for the classes to be multiplied and stacks the boxed links on top of each other, with the right boxed link at the bottom. Then one erases the line in the middle which results in a boxed linked, {\it i.e.} in an new equivalence class. Note that the unit of the braid group is represented by $N$ parallel boxed lines.

\vspace{0.2cm}

As it is apparent from the above diagrams, the braid operations can be generated from a set of basic moves, called the generators. For $\Bb_4$, there are 3 generators and they look as below:
\begin{equation}\label{Eq:B4Gen}
\beta_1=\begin{matrix} 
\begin{tikzpicture}
\braid[number of strands=4,width=0.4cm,height=0.4cm,border height=0.1cm] a_1;
\draw (0.2,-0.6) rectangle (1.8,0);
\end{tikzpicture}
\end{matrix}
\
\begin{matrix} 
, \quad
 \end{matrix}
 \
\beta_2 = \begin{matrix} 
\begin{tikzpicture}
\braid[number of strands=4,width=0.4cm,height=0.4cm,border height=0.1cm] a_2;
\draw (0.2,-0.6) rectangle (1.8,0);
\end{tikzpicture} 
\end{matrix}
\
\begin{matrix} 
, \quad
\end{matrix}
 \
\beta_3 = \begin{matrix} 
\begin{tikzpicture}
\braid[number of strands=4,width=0.4cm,height=0.4cm,border height=0.1cm] a_3;
\draw (0.2,-0.6) rectangle (1.8,0);
\end{tikzpicture} 
\end{matrix} 
\end{equation}
The inverse of these generators are:
\begin{equation}\label{Eq:B4GenInv}
\beta_1^{-1}=\begin{matrix} 
\begin{tikzpicture}
\braid[number of strands=4,width=0.4cm,height=0.4cm,border height=0.1cm] a_1^{-1};
\draw (0.2,-0.6) rectangle (1.8,0);
\end{tikzpicture}
\end{matrix}
\
\begin{matrix} 
, \quad
 \end{matrix}
 \
\beta_2^{-1} = \begin{matrix} 
\begin{tikzpicture}
\braid[number of strands=4,width=0.4cm,height=0.4cm,border height=0.1cm] a_2^{-1};
\draw (0.2,-0.6) rectangle (1.8,0);
\end{tikzpicture} 
\end{matrix}
\
\begin{matrix} 
, \quad
\end{matrix}
 \
\beta_3^{-1} = \begin{matrix} 
\begin{tikzpicture}
\braid[number of strands=4,width=0.4cm,height=0.4cm,border height=0.1cm] a_3^{-1};
\draw (0.2,-0.6) rectangle (1.8,0);
\end{tikzpicture} 
\end{matrix} 
\end{equation}
Indeed, a quick check shows:
\begin{equation}\label{Eq:GenCheck}
\begin{matrix} 
\begin{tikzpicture}
\braid[number of strands=4,width=0.4cm,height=0.4cm,border height=0.1cm] a_1;
\draw (0.2,-0.6) rectangle (1.8,0);
\end{tikzpicture}
\end{matrix}
\
\begin{matrix} 
\cdot
 \end{matrix}
 \
 \begin{matrix} 
\begin{tikzpicture}
\braid[number of strands=4,width=0.4cm,height=0.4cm,border height=0.1cm] a_1^{-1};
\draw (0.2,-0.6) rectangle (1.8,0);
\end{tikzpicture}
\end{matrix}
\
\begin{matrix} 
=
 \end{matrix}
 \
  \begin{matrix} 
\begin{tikzpicture}
\braid[number of strands=4,width=0.4cm,height=0.4cm,border height=0.1cm] a_1 a_1^{-1};
\draw (0.2,-1) rectangle (1.8,0);
\end{tikzpicture}
\end{matrix}
\end{equation}
In general, $\Bb_N$ has $N-1$ generators and in terms of these generators, which will be denoted by $\sigma_1$, \ldots, $\sigma_{N-1}$, the braid group can be defined as the free group generated by the $\sigma$'s modulo the following relations:
\begin{equation}\label{Eq:ArtinRelations}
\begin{array}{ll}
& \sigma_j \sigma_{j+1} \sigma_j = \sigma_{j+1} \sigma_j \sigma_{j+1}, \quad  j=1,\ldots,N-2, \\
& \sigma_j \sigma_k = \sigma_k \sigma_j, \quad j,k=1,\ldots,N-1, \quad |j-k| \geq 2.
\end{array}
\end{equation}
This is known as the Artin presentation of the braid group \cite{ArtinAM1947}. The first set of relations translates the equivalence of \eqref{Eq:121} and \eqref{Eq:212}. Let us point out that, although the braid group has finite generators, the group itself is of infinite cardinality. Indeed, in the Artin presentation, an element of $\Bb_N$ takes the form:
\begin{equation}
\big [\sigma_{i_1}^{p_1} \cdot \sigma_{i_2}^{p_2} \ldots \sigma_{i_K}^{p_K} \big ], \quad K \in \NM_+,
\end{equation}
where the indices $i_\alpha$ are drawn from $\{1,\ldots,N-1\}$ and the powers $p_\alpha$ can be any natural number. The convention is that $\sigma_{i_\alpha}^0 = 1$, the identity of the braid group. The square brackets are there to remind us that these monoidals are defined only up to the equivalence relations \eqref{Eq:ArtinRelations}. In other words, that we are dealing with an equivalence class.

\subsection{Planar Diagrams} A planar $N$-diagram consists of:
\begin{itemize}
\item A set of $N$ equidistant points located at the bottom of a box;
\item A set of $N$ equidistant points located at the top of a box;
\item $N$ smooth planar strands emerging from a point at the bottom or at the top of a box, and diving in a point at the top as well as at the bottom of same box, always normally to the edges of the box. 
\item Arbitrary number of strands that close into themselves and are suspended inside the box;
\end{itemize}
None of the strands are allowed to intersect the other and the strands are defined up to deformations allowed by topology. Below are examples of such diagrams:
\begin{equation}\label{Eq:PlanarD1}
\begin{matrix}
\begin{tikzpicture}
    \draw (0,0) rectangle (3.5,1);
    \draw (0.5,0)  .. controls +(0,0.3) and +(0,0.3) .. (1,0);
    \draw (0.5,1)  .. controls +(0,-0.3) and +(0,0.3) .. (1.5,0);
    \draw (1,1)  .. controls +(0,-0.3) and +(0,-0.3) .. (1.5,1);
    \draw (2,0)  -- (2,1);
    \draw (2.5,1)  .. controls +(0,-0.3) and +(0,-0.3) .. (3,1);
    \draw (2.5,0)  .. controls +(0,0.3) and +(0,0.3) .. (3,0);
    \draw (1.7,0.5) circle (0.15);
    \draw (2.3,0.5) circle (0.15);
  \end{tikzpicture}
  \end{matrix} \ 
  \begin{matrix}
  ,
  \end{matrix} \quad 
  \begin{matrix}
\begin{tikzpicture}
    \draw (0,0) rectangle (3.5,1);
    \draw (0.5,1) .. controls +(0,-0.3) and +(0,-0.3) .. (1,1);
    \draw (0.5,0)  .. controls +(0,0.3) and +(0,-0.3) .. (1.5,1);
    \draw (1,0)  .. controls +(0,0.3) and +(0,-0.3) .. (2,1);
    \draw (1.5,0)  .. controls +(0,0.3) and +(0,0.3) .. (2,0);
    \draw (2.5,1)  .. controls +(0,-0.3) and +(0,-0.3) .. (3,1);
    \draw (2.5,0)  .. controls +(0,0.3) and +(0,0.3) .. (3,0);
  \end{tikzpicture}
  \end{matrix}
\end{equation}
Note that, necessarily, every single point at the top or at the bottom of the box are paired by a strand with one and only one other such point. Again, we should think of a specific planar diagram as a representative of all its equivalent diagrams and by a boxed planar diagram as the equivalence class of that particular diagram. 

\vspace{0.2cm}

Planar diagrams can be composed in a similar manner as the link diagrams. For example, the composition of the diagrams shown in \eqref{Eq:PlanarD1} gives:
\begin{equation}\label{Eq:PlanarD1}
\begin{matrix}
\begin{tikzpicture}
    \draw (0,0) rectangle (3.5,1);
    \draw (0.5,0)  .. controls +(0,0.3) and +(0,0.3) .. (1,0);
    \draw (0.5,1)  .. controls +(0,-0.3) and +(0,0.3) .. (1.5,0);
    \draw (1,1)  .. controls +(0,-0.3) and +(0,-0.3) .. (1.5,1);
    \draw (2,0)  -- (2,1);
    \draw (2.5,1)  .. controls +(0,-0.3) and +(0,-0.3) .. (3,1);
    \draw (2.5,0)  .. controls +(0,0.3) and +(0,0.3) .. (3,0);
    \draw (1.7,0.5) circle (0.15);
    \draw (2.3,0.5) circle (0.15);
  \end{tikzpicture}
  \end{matrix} \
  \begin{matrix}
  \cdot
  \end{matrix}\
  \begin{matrix}
\begin{tikzpicture}
    \draw (0,0) rectangle (3.5,1);
    \draw (0.5,1) .. controls +(0,-0.3) and +(0,-0.3) .. (1,1);
    \draw (0.5,0)  .. controls +(0,0.3) and +(0,-0.3) .. (1.5,1);
    \draw (1,0)  .. controls +(0,0.3) and +(0,-0.3) .. (2,1);
    \draw (1.5,0)  .. controls +(0,0.3) and +(0,0.3) .. (2,0);
    \draw (2.5,1)  .. controls +(0,-0.3) and +(0,-0.3) .. (3,1);
    \draw (2.5,0)  .. controls +(0,0.3) and +(0,0.3) .. (3,0);
  \end{tikzpicture}
  \end{matrix}\
  \begin{matrix}
  =
  \end{matrix}\
  \begin{matrix}
\begin{tikzpicture}
    \draw (0,0) rectangle (3.5,2);
    \draw (0.5,1) .. controls +(0,-0.3) and +(0,-0.3) .. (1,1);
    \draw (0.5,0)  .. controls +(0,0.3) and +(0,-0.3) .. (1.5,1);
    \draw (1,0)  .. controls +(0,0.3) and +(0,-0.3) .. (2,1);
    \draw (1.5,0)  .. controls +(0,0.3) and +(0,0.3) .. (2,0);
    \draw (2.5,1)  .. controls +(0,-0.3) and +(0,-0.3) .. (3,1);
    \draw (2.5,0)  .. controls +(0,0.3) and +(0,0.3) .. (3,0);
    \draw (0.5,1)  .. controls +(0,0.3) and +(0,0.3) .. (1,1);
    \draw (0.5,2)  .. controls +(0,-0.3) and +(0,0.3) .. (1.5,1);
    \draw (1,2)  .. controls +(0,-0.3) and +(0,-0.3) .. (1.5,2);
    \draw (2,1)  -- (2,2);
    \draw (2.5,2)  .. controls +(0,-0.3) and +(0,-0.3) .. (3,2);
    \draw (2.5,1)  .. controls +(0,0.3) and +(0,0.3) .. (3,1);
    \draw (1.7,1.5) circle (0.15);
    \draw (2.3,1.5) circle (0.15);
  \end{tikzpicture}
  \end{matrix}
\end{equation}
The composition with $N$ straight strands does not modify the class of a planar diagram, hence it serves again as the unit. However, not any planar diagram has an inverse, hence the set of planar diagrams and the composition rule forms only a monoid. 

\vspace{0.2cm}

Any planar diagram which contains closed loops can be written as the product of $N$ straight strands plus all the closed loops and a diagram with no closed loops. For example:
\begin{equation}\label{Eq:DElelemnt1}
\begin{matrix}
\begin{tikzpicture}
    \draw (0,0) rectangle (3.5,1);
    \draw (0.5,0)  .. controls +(0,0.3) and +(0,0.3) .. (1,0);
    \draw (0.5,1)  .. controls +(0,-0.3) and +(0,0.3) .. (1.5,0);
    \draw (1,1)  .. controls +(0,-0.3) and +(0,-0.3) .. (1.5,1);
    \draw (2,0)  -- (2,1);
    \draw (2.5,1)  .. controls +(0,-0.3) and +(0,-0.3) .. (3,1);
    \draw (2.5,0)  .. controls +(0,0.3) and +(0,0.3) .. (3,0);
    \draw (1.7,0.5) circle (0.15);
    \draw (2.3,0.5) circle (0.15);
  \end{tikzpicture}
  \end{matrix} \
  \begin{matrix}
  =
  \end{matrix}\
  \begin{matrix}
\begin{tikzpicture}
    \draw (0.25,0.3) circle (0.15);
    \draw (0.25,0.7) circle (0.15);
    \draw (0,0) rectangle (3.5,1);
    \draw (0.5,0)  -- (0.5,1);
    \draw (1,0)  -- (1,1);
    \draw (1.5,0)  -- (1.5,1);
    \draw (2,0)  -- (2,1);
    \draw (2.5,0) -- (2.5,1);
    \draw (3,0)  -- (3,1);
  \end{tikzpicture} 
  \end{matrix}\
  \begin{matrix}
  \cdot
  \end{matrix} \
  \begin{matrix}
\begin{tikzpicture}
    \draw (0,0) rectangle (3.5,1);
    \draw (0.5,0)  .. controls +(0,0.3) and +(0,0.3) .. (1,0);
    \draw (0.5,1)  .. controls +(0,-0.3) and +(0,0.3) .. (1.5,0);
    \draw (1,1)  .. controls +(0,-0.3) and +(0,-0.3) .. (1.5,1);
    \draw (2,0)  -- (2,1);
    \draw (2.5,1)  .. controls +(0,-0.3) and +(0,-0.3) .. (3,1);
    \draw (2.5,0)  .. controls +(0,0.3) and +(0,0.3) .. (3,0);
  \end{tikzpicture}
  \end{matrix}
  \end{equation}
Furthermore, the element:
\begin{equation}\label{Eq:DElement}
\begin{matrix}
d
\end{matrix}\
\begin{matrix}
=
\end{matrix}\
  \begin{matrix}
\begin{tikzpicture}
    \draw (0.25,0.5) circle (0.15);
    \draw (0,0) rectangle (3.5,1);
    \draw (0.5,0)  -- (0.5,1);
    \draw (1,0)  -- (1,1);
    \draw (1.5,0)  -- (1.5,1);
    \draw (2,0)  -- (2,1);
    \draw (2.5,0) -- (2.5,1);
    \draw (3,0)  -- (3,1);
  \end{tikzpicture} 
  \end{matrix}
  \end{equation}
  commutes with all the elements of the monoid and \eqref{Eq:DElelemnt1} can be written as:
  \begin{equation}\label{Eq:DElelemnt2}
\begin{matrix}
\begin{tikzpicture}
    \draw (0,0) rectangle (3.5,1);
    \draw (0.5,0)  .. controls +(0,0.3) and +(0,0.3) .. (1,0);
    \draw (0.5,1)  .. controls +(0,-0.3) and +(0,0.3) .. (1.5,0);
    \draw (1,1)  .. controls +(0,-0.3) and +(0,-0.3) .. (1.5,1);
    \draw (2,0)  -- (2,1);
    \draw (2.5,1)  .. controls +(0,-0.3) and +(0,-0.3) .. (3,1);
    \draw (2.5,0)  .. controls +(0,0.3) and +(0,0.3) .. (3,0);
    \draw (1.7,0.5) circle (0.15);
    \draw (2.3,0.5) circle (0.15);
  \end{tikzpicture}
  \end{matrix} \
  \begin{matrix}
  =
  \end{matrix}\
  \begin{matrix}
d^2 
  \end{matrix}\
  \begin{matrix}
\begin{tikzpicture}
    \draw (0,0) rectangle (3.5,1);
    \draw (0.5,0)  .. controls +(0,0.3) and +(0,0.3) .. (1,0);
    \draw (0.5,1)  .. controls +(0,-0.3) and +(0,0.3) .. (1.5,0);
    \draw (1,1)  .. controls +(0,-0.3) and +(0,-0.3) .. (1.5,1);
    \draw (2,0)  -- (2,1);
    \draw (2.5,1)  .. controls +(0,-0.3) and +(0,-0.3) .. (3,1);
    \draw (2.5,0)  .. controls +(0,0.3) and +(0,0.3) .. (3,0);
  \end{tikzpicture}
  \end{matrix}
  \end{equation}
  As such, any planar diagram can be expressed as the product of $d$ to some power and a planar diagram with no closed loops inside it.
  
 \vspace{0.2cm}
 

We call the $n$ dimensional vector space of planar diagrams $TL_n (d)$. For a given dimension $n$, there are $\frac{1}{n+1}$ $2n \choose n$ distinct planar diagrams which form a basis of $TL_n (d)$

\section{The code} 
\subsection{Generating Planar Diagrams}
We want to generate the space of all possible planar diagrams (PD) using a simple code. Since planar diagrams are primarily pairwise relationships between two nodes, along with some additional conditions which restrict how connections can be made, we consider a matrix $P$  of size $2N$ by $2N$ where $N$ is the number of bottom nodes in the planar diagram (making the total number of nodes $2N$)  with $P_{ij} = 1$ if node $i$ is connected to node $j$ on the planar diagram and $0$ otherwise. If one were to consider all pairwise connections between two nodes in a $2N$ node system, the number of possible diagrams to consider would be of size $\frac{(2N)!}{N!}$. Due the the inherent restrictions governing which connections are possible in planar diagrams, we in fact only have $\frac{(2N)!}{(N+1)!(N!)^2}$ many unique diagrams. It is therefore importent to consider how those rescrictions translate to a smaller possibility space. 

\begin{figure}[h]
	\includegraphics[height=3cm, width=3cm]{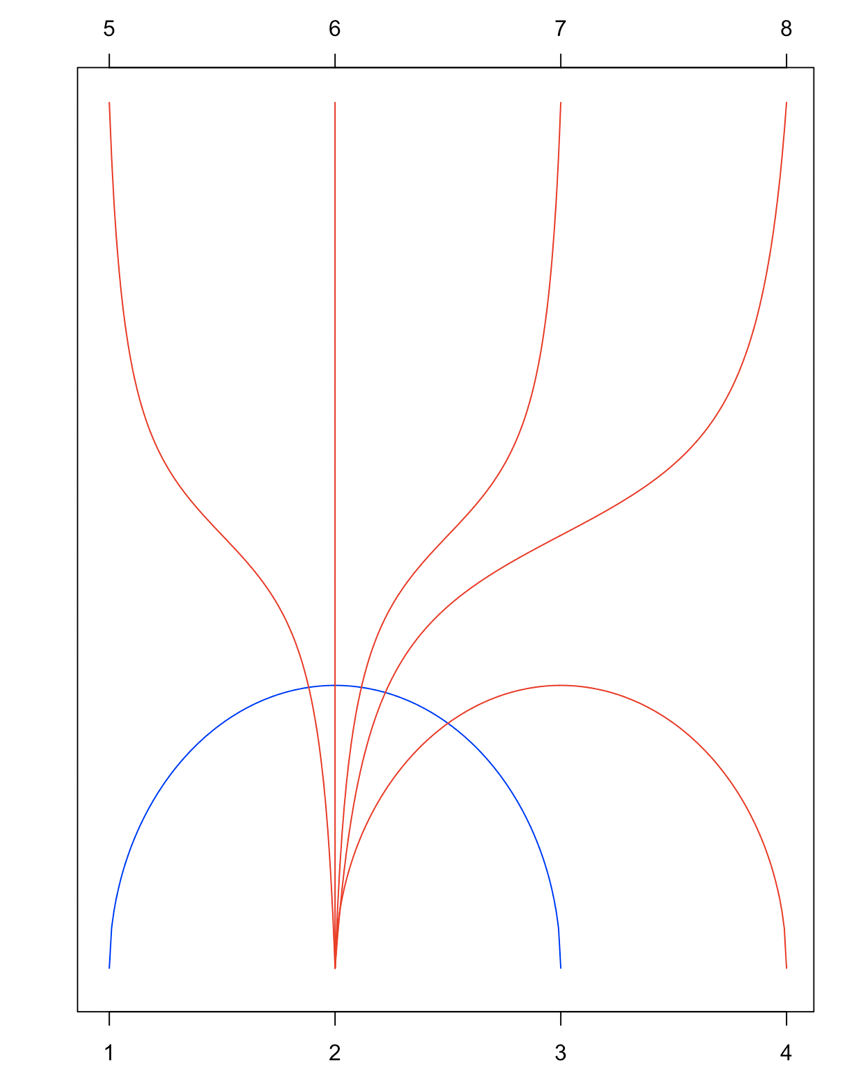}
	\includegraphics[height=3cm, width=3cm]{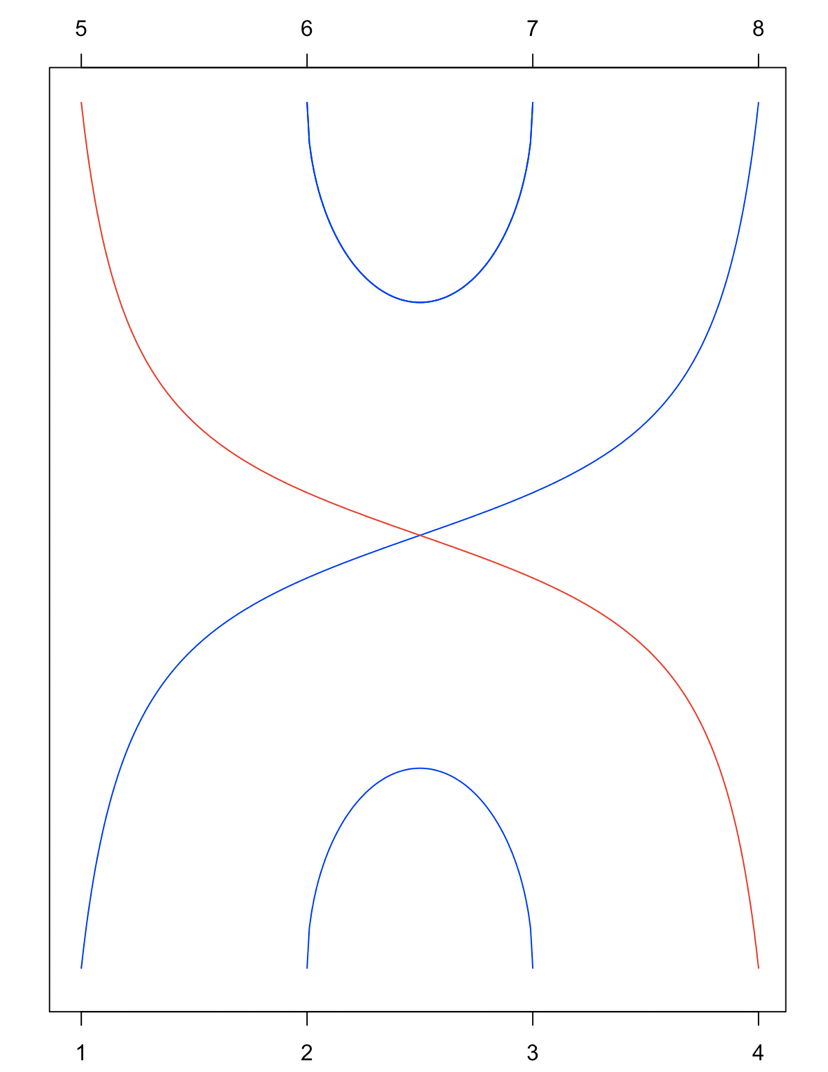}
	\caption{On the left we have that no matter the available options, node $2$ cannot connect without an edge intersection, and on the right we see that whichever nodes are left out from the remaining top-top/bottom-bottom connections, are forced to have their edges intersect.}
	\label{PossibleConnections}
\end{figure}

From the start it is clear that if a planar diagram is to have no edge intersections it must be the case that not every node can connect to any other node. Take for example an edgeless PD with $4$ nodes on the bottom and $4$ on the top, with the bottom nodes labeled $1$ through $4$, and the top labeled $5$ through $8$ (from left to right). While no edges have yet been placed, it remains the case that node $1$ is only connectable with nodes $\{2,4,5,7\}$. For if it were to connect with node $3$, for example, node $2$ would be forced to have its edge intersect with the $1<->3$ edge or to remain connectionless, both of which options violate the definition of a PD.
Further, node $1$ cannot connect with node $8$ since that would prevent any other possible top-bottom connections, leaving an odd number of nodes both on top and on bottom, implying at least two nodes which can never be connected. See Figure \ref{PossibleConnections}.
Thus the set of nodes to which node $i$ can connect to even on an edgeless PD is a strict subset of the full $2N$ node set, containing further a dependence on $N$ itself. In particular, let $\gamma(i,j)$ be the connectability of nodes $i$ and $j$, with $$  \begin{cases}
\gamma(i,j) = 1 & \text{if nodes $i$ and $j$ are connectable }\\ 
\gamma(i,j) = 0& \text{if nodes $i$ and $j$ are not connectable}
\end{cases} $$ Then for a given $N$, if $N$ is odd, we have $	\gamma(i,j) = 1$ if $|i-j|$ is odd, and if $N$ is even then, if $i$ and $j$ are either both on top or on bottom $	\gamma(i,j) = 1$ when $|i-j|$ is odd, otherwise $\gamma(i,j) = 1$ when $|i-j|$ is even. 
\begin{code}{\rm Below are code lines which performs the binary connectablility of any two nodes  $n \in \{0,1,\ldots, 2N\}$.
		\begin{equation}\label{Initial_connectablility}
		\boxed{
			\begin{array}{l}
			\ \  \mathrm{\gamma} = \mathrm{matrix_{2N,2N}} \\
			\ \ \mathrm{for} \ i \ \in \mathrm{1:2N} \\
	        \ \	\ \ \mathrm{for} \ j \ \in \mathrm{1:2N} \\
			\ \ \ \ \ \ \mathrm{if \ i,j \in 1:N} \ \mathrm{or \ i,j \in (N+1):2N} \\
			\ \ \ \ \ \ \ \ \mathrm{if \ {abs(i-j) (mod 2) }= 0 } \\
			\ \ \ \ \ \ \ \ \ \ \mathrm{\gamma_{i,j} = 0 } \\
			\ \ \ \ \ \ \ \ \mathrm{else } \\
			\ \ \ \ \ \ \ \ \ \ \mathrm{\gamma_{i,j} = 1 } \\
	     	\ \ \ \ \ \ \mathrm{else} \\
	     	\ \ \ \ \ \ \ \ \mathrm{if \ {abs(i-j) (mod 2)}= 0 } \\
	     	\ \ \ \ \ \ \ \ \ \ \mathrm{\gamma_{i,j} = (N+1)  (mod 2)} \\
	     	\ \ \ \ \ \ \ \ \mathrm{else } \\
	     	\ \ \ \ \ \ \ \ \ \ \mathrm{\gamma_{i,j} = N  (mod 2)} \\
	     	\mathrm{return} \gamma
			\end{array}
		}
		\end{equation}
	} 
\end{code}

\begin{figure}[h]
	\includegraphics[height=3cm, width=3cm]{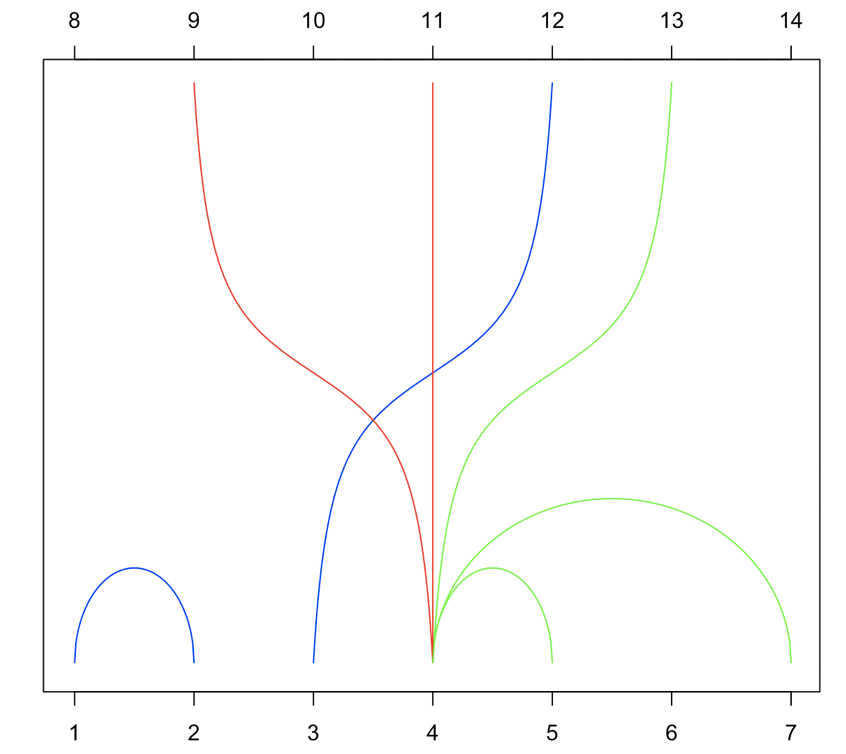}
	\includegraphics[height=3cm, width=3cm]{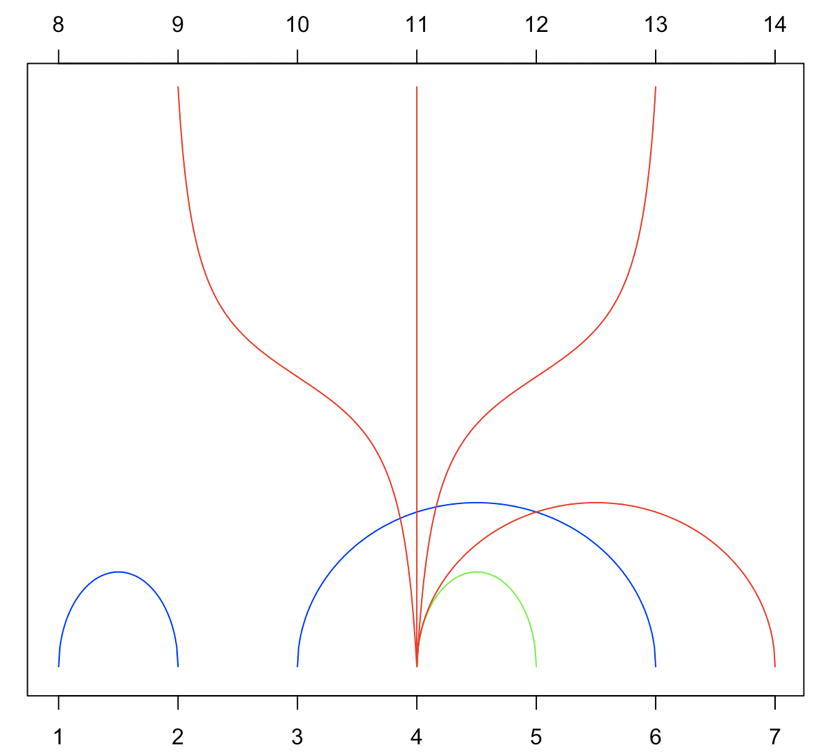}
	\caption{On the left we see that early bottom-top connections (seen in blue) can restrict a later node's ability to have bottom top connections, in particular when the node with which it connects is less than the maximum top node already connected (seen in red). They however do not restrict any other possible connections (seen in green). On the right we can see that any bottom-bottom connections surrounding $i$ will result in $i$ being able to connect only with those within the same bubble.}
	\label{EdgeIntersections}
\end{figure}

Now that it is clear which nodes can connect to which, even before any edges have been placed (in anticipation of future edge intersections), we turn to considering the particular effects which previously placed edges have on future possibilities. Namely, let $PD^i$ be a planer diagram with fixed edges connecting all nodes up to node $i \leq N$ (not including $i$) and let $x  = \{1,2,3, ..., i-1\}$ be the set of connected nodes. Finally, let $P^i$ be the $2N$ by $2N$ matrix associated with $PD^i$, where $P^i_{j,k} = 1$ if node $j$ is connected with node $k$ (by definition $P^i_{j,k} = 1$ implies either $j \in x$ or $k \in x$, and $P^i_{j,k} = 0$ if $j \& k \geq i$  ). It is clear that node $i$ cannot connect to any $j \in x$ even if $\gamma(i,j) = 1$, for they have already been connected by assumption. Further, if node $1$ connects to node $i+N$ (which is opposite node $i$ on the upper level of the PD), then those nodes which are greater than $N$ but less than $N+i$ are also no longer connectable with $i$. Indeed if any of the nodes in $x$ connect to node $N+i$, we would have that $i$ is no longer connectable with nodes greater than $N$ but less than $N+i$. We could generalize even further by recognizing that any node $j > N$ will not be connectable with $i$ if there is a node in $x$ which connects to a node $k > j$, see Figure \ref{EdgeIntersections}.
 Lastly we have the edge restrictions which arise from having a node $j \in x$ connect with a node $k$ with $i \leq k \leq N$  (bottom-bottom connections which surround the node $i$). Similar to before, these edge placements will make it such that $i$ can only connect to nodes which are less than the minimum such $k$ and greater than $i$, see Figure \ref{EdgeIntersections}.
 \begin{code}{\rm Below are code lines which performs the removal of connectablility between node $i$  in $PD^i$ defined above and $j\in  \{1,2,3, \dots , 2N\} $.
 		\begin{equation}\label{edge_removal}
 		\boxed{
 			\begin{array}{l}
 			\mathrm{input} \ PD^i \\
 			x = \{1,2,3, \dots , i-1\}\\
 			y = which(P^i_{x,:} == 1) \\
 			\gamma_{i,j} = 0 \ \ \mathrm{for} \ j \in (x \cup y)\\
 			y_{max} = max(y)\\
 			\mathrm{if} \ y_{max} > N \\
 			\ \ 	\gamma_{i,j} = 0 \ \ \mathrm{for} \ N+1 \leq j \leq y_{max}\\
 			y_{min} = min(y)\\
 			\mathrm{if} \ y_{min} \leq N \\
 			\ \ 	\gamma_{i,j} = 0 \ \ \mathrm{for} \ j \geq y_{min} \\
 			\mathrm{return} \gamma
 			\end{array}
 		}
 		\end{equation}
 	} 
 \end{code}
  It is worthwhile to note that though in principle we should consider bottom-top  connections to the right of node $i$ as we did for those to the left of node $i$, if one places edges sequentially (from $i= 1$ to $i = 2N$) then future connections need not be considered.

\begin{figure}[h]
	\includegraphics[height=1cm, width=3cm]{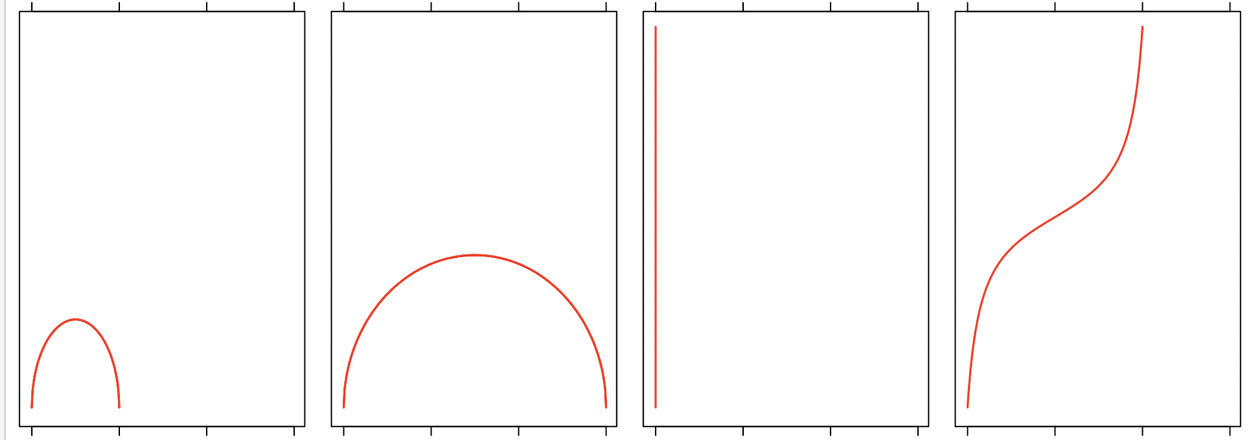}\\
	\includegraphics[height=1cm, width=5cm]{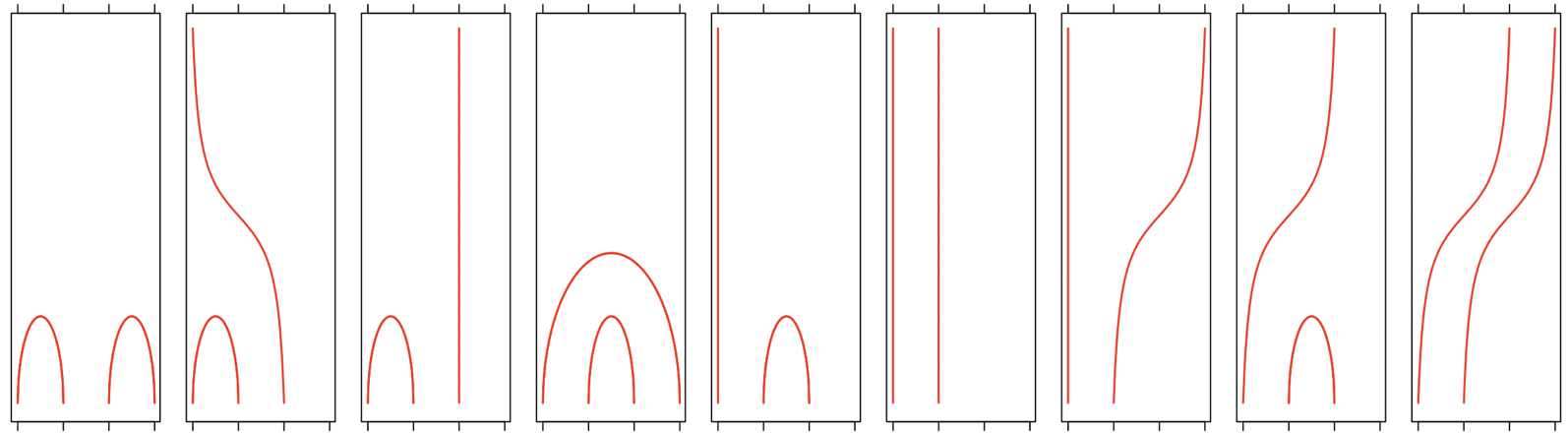}\\
	\includegraphics[height=1cm, width=8.5cm]{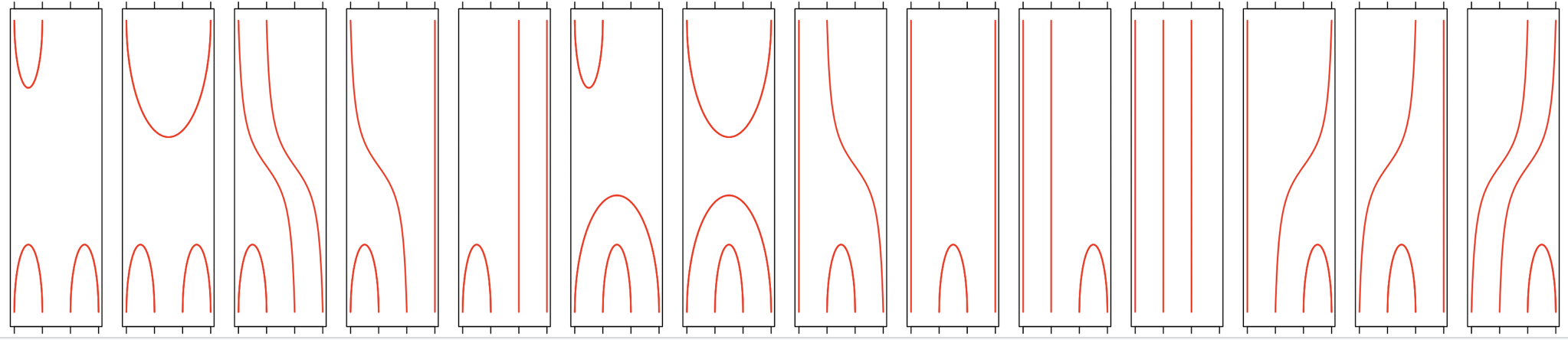}\\
	\includegraphics[height=1  cm, width=9.5cm]{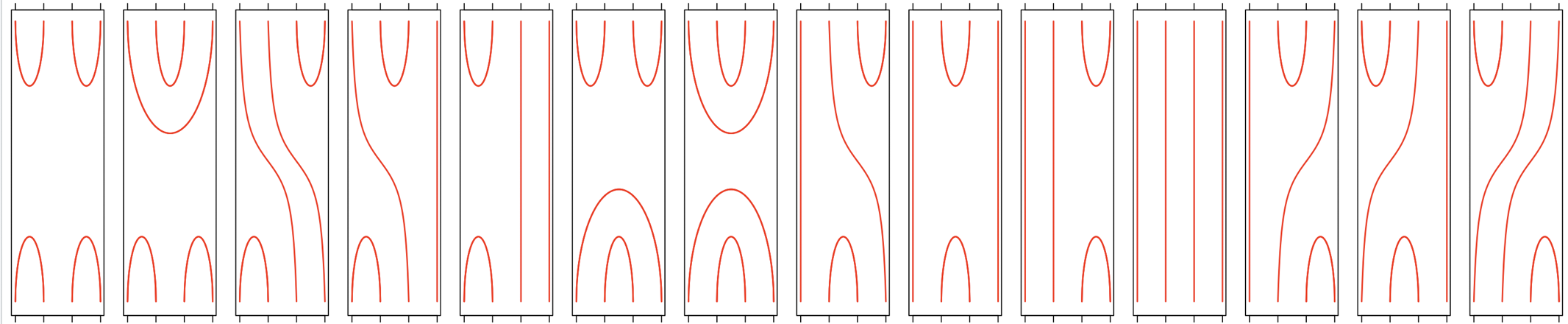}\\

	\caption{From top to bottom we can see the iterative generation of the planer diagrams with the first row being all possible first connections of node $1$ and the second row being all possible connections of node $2$ for each of the possibilities in row $1$ and so forth for rows $3$ and $4$. }
	\label{ItterativeGenerating}
\end{figure}

\begin{figure}
	\includegraphics[height=10cm, width=10cm]{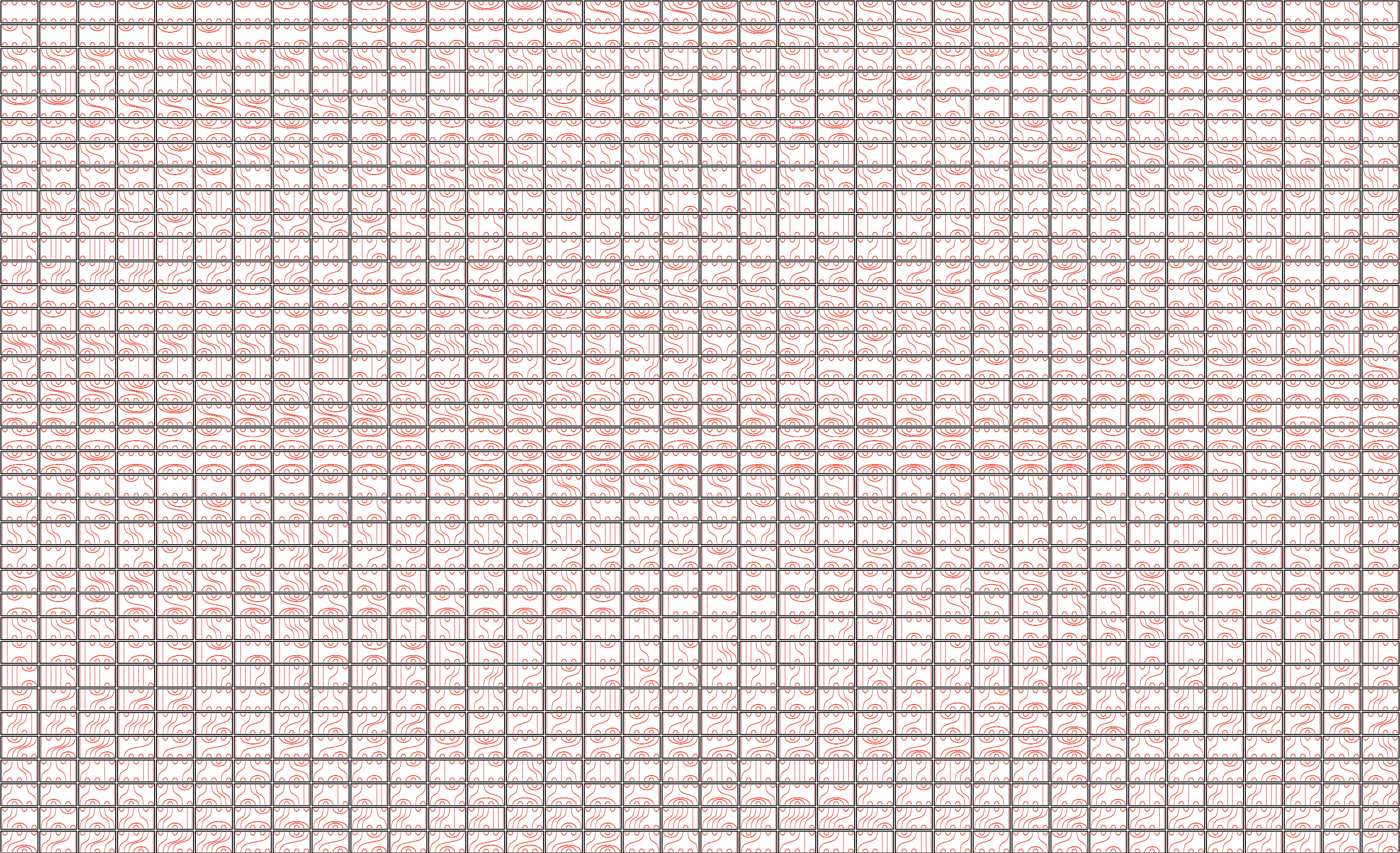}
	
	\caption{Seen above is the set of all diagrams for N = 8.}
	\label{N=8Diagrams}
\end{figure}

With these restrictions in mind, we can now begin generating all of the possible PDs sequentially. Beginning with node $1$, we have $\kappa(1)$,  the set of all possible PDs which have nothing but the first node connected. Form the above discussion it is clear that the cardinality of $\kappa(1)$ is $N$ (since there are $N$ possible nodes to which node $1$ can connect). If $x \in \kappa(1)$ is one particular PD, then let $\phi(x)$ be all possible planer diagrams which can be generated from $x$ by adding only one edge (since the number of next possible edges depends on $x$, so will the cardinality of $\phi(x)$ depend on $x$). We can now iteratively span the possible PDs via,
$$\kappa(i)  = \{ y = \phi(x) | \forall x \in \kappa(i-1)\}$$
with $\kappa(N)$ being the complete set of planer diagrams with $N$ bottom nodes see Figures \ref{ItterativeGenerating} and \ref{N=8Diagrams}. The code was verified via the equation  $\frac{1}{n+1}$ $2n \choose n$ and by counting the number of PDs obtained for $N \leq 8$.  

\begin{code}{\rm Below are code lines which sequentially span the set of possible PDs for dim = N.
		\begin{equation}\label{Spanning}
		\boxed{
			\begin{array}{l}
			\gamma = \mathrm{Code_{\ref{Initial_connectablility}}}\\
			x  = which(\gamma[1,:] == 1)\\
			\mathrm{counte}r = 1\\
			\mathrm{for} \ i \ \in \  x \\
			\ \ \kappa(1)(\mathrm{counter})[1,i] =1 \\
			\ \ \mathrm{counter} = \mathrm{counter} +1\\
			\mathrm{counter_1} = 1\\
			\mathrm{while \ counter_1 < N} \\
			\ \ \mathrm{counter_2} = 1\\
			\ \  \mathrm{counter_3} = 1\\
			\ \ \mathrm{while \ counter_2 <\mathrm{length(\kappa(\mathrm{counter_1}))} } \\
		 	 \ \  \ \ \gamma = \mathrm{Code}_{\ref{edge_removal}}(\kappa(\mathrm{counter_1})(\mathrm{counter_2}))\\
			\ \ \ \ x = which (\gamma_{\mathrm{counter_1+1}, :} ==1 )\\
		
		    \ \	\ \ \mathrm{for} \ j \ \in \ x\\
			\ \ \ \ \ \ \kappa(\mathrm{counter_1+1})(\mathrm{counter_3})[\mathrm{counter_1+1}, j] =1 \\
			\ \ \ \ \ \ \mathrm{counter_3} =\mathrm{counter_3}+ 1\\
			 \ \ \ \ \mathrm{counter_2 = counter_2 +1}\\
			 \ \ \mathrm{counter_1} = \mathrm{counter_1 +1}
			
			\end{array}
		}
		\end{equation}
	} 
\end{code}

\newpage
\subsection{The representation of Planer Diagram multiplication} 
\begin{figure}
\includegraphics[height=2cm, width=2cm]{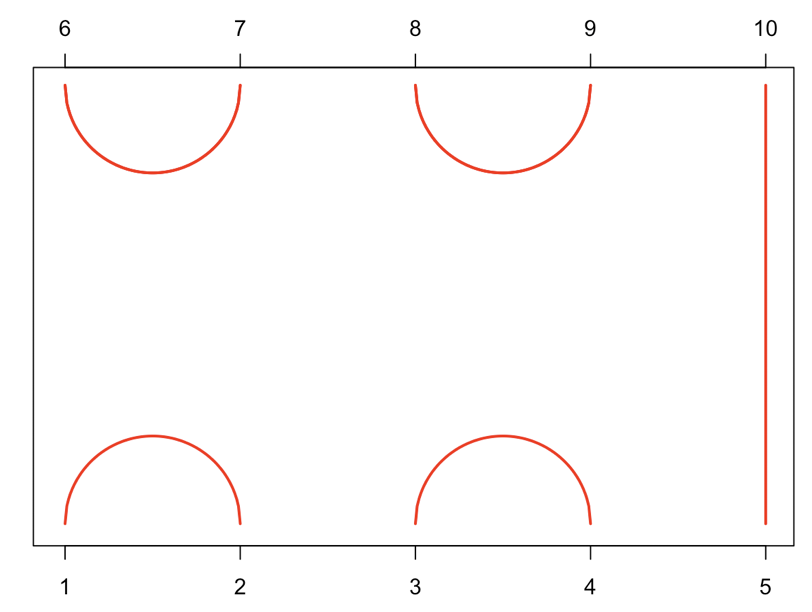}
\includegraphics[height=2cm, width=2cm]{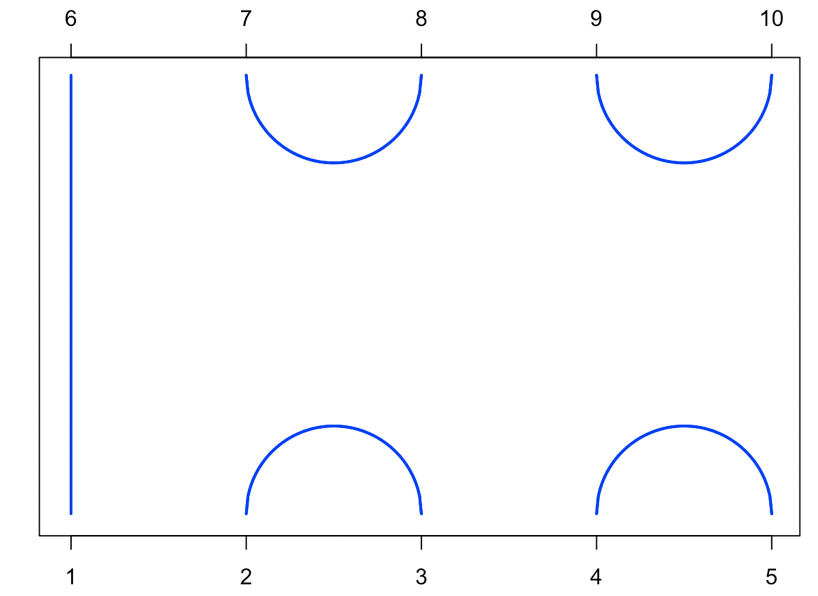} \ \
\includegraphics[height=2cm, width=2cm]{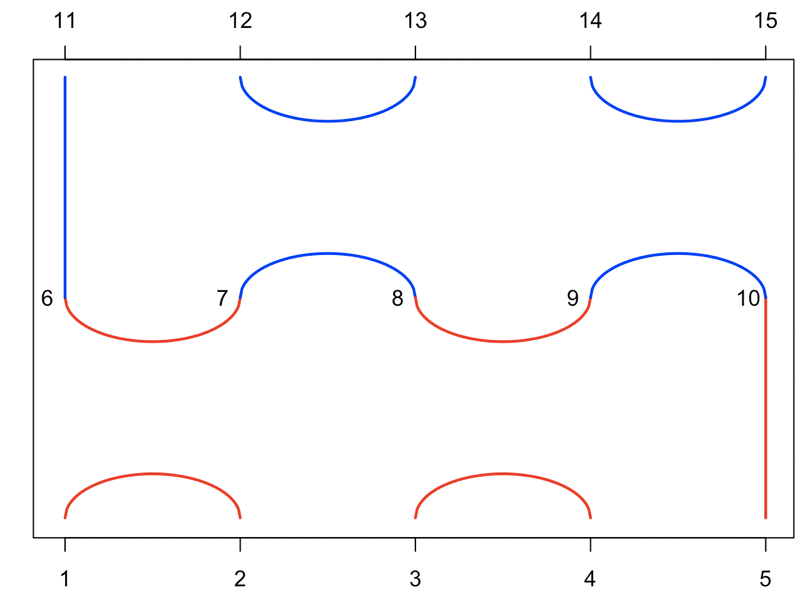}
\includegraphics[height=2cm, width=2cm]{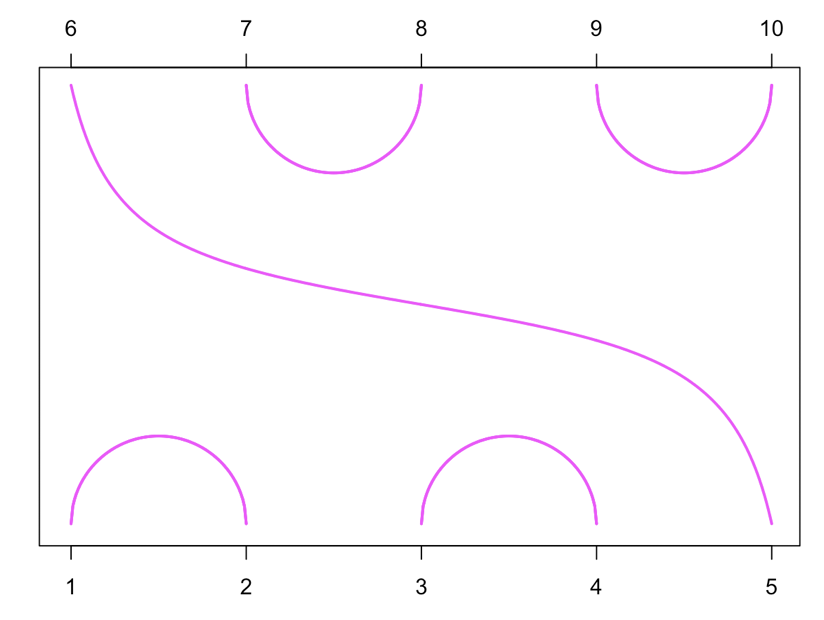} \ \
\includegraphics[height=2cm, width=2cm]{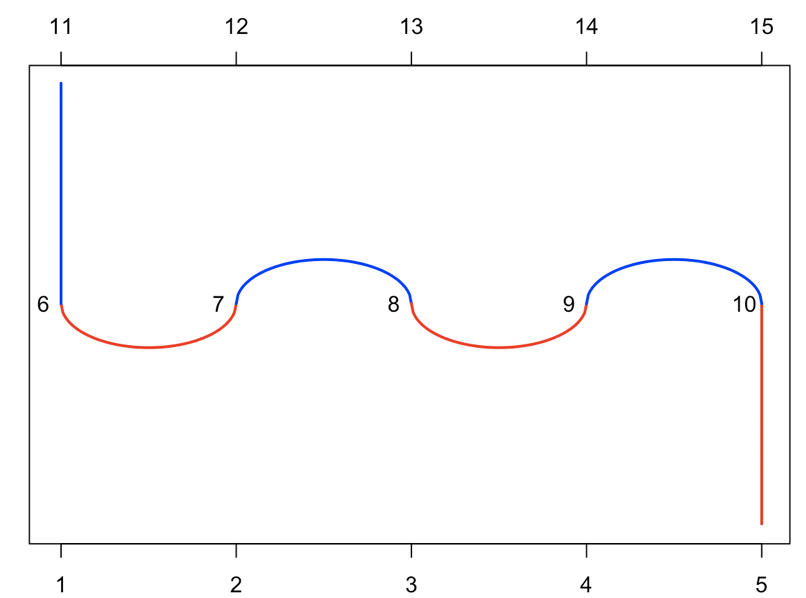}

\caption{Seen above is an example of a planer diagram multiplication. Top row shows $D^1$(red), $D^2$(blue), $D^1D^2$ (both mid operation (red/blue) and final diagram (pink)) and the path connecting nodes 5 and 11, seen from left to right. Note, nodes 5 and 6 are connected in the product diagram  by vitue of the fact that there is a set $c = \{10,9,8,7,6\}$ seen in the right most graph satisfying condition (3) mentioned.}
\label{ExampleProduct}
\end{figure}

Now that we can both ariticulate planer diagrams using matricies and generate the space for any given dimention, we would like to consider the function which produces the product diagram for any two diagrams. As mentioned in the previus section, the product diagram of two diagrams is the sequential stacking of the second onto the first. In particular, we can consider all pairwise connections in the product diagram:($\{(a,b)| a,b \leq N \}$, $\{(a,b)| a\leq N < b \}$, $\{(a,b)| a,b > N \}$) and determine whether the pair $(a,b)$ is connected, given the two input diagrams. \\
Let $D^1, D^2$ be two planer diagrams, with $A^1,A^2$ thier respective matrix representations, and let $D^3 = D^1  D^2$ with $A^3$, its matrix representation. If $(a,b)$ are a pair of nodes in $D^3$ then $A^3_{a,b} = 1$ if
\begin{enumerate}
	\item $A^1_{a,b} = 1$ when $a,b \leq N$
	\item $A^2_{a,b} = 1$ when $a,b > N$
	\item If there exists a sequence $(c_1,c_2,c_3...c_n)| c_i \in (N+1,2N) $ satisfying 
	\begin{itemize}
		\item $A^1_{a,c_1} = 1$
		\item $A^1_{c_{i+1},c_i}+A^2_{c_{i+1}-N,c_i-N} > 0 \ \ \ \ \forall i \in [1,n]$
		\item $A^2_{c_n-N,b} = 1$
	\end{itemize}
	 when $a \leq N< b$
\end{enumerate}
In other words, two nodes are connected in the product diagram if there is a path connecting them when the second diagram is placed ontop of the first see Figure \ref{ExampleProduct}. This can also be verified by evaluating the product series $P = \prod_{i=1}^{n-1} (A^1_{c_{i+1},c_i}+A^2_{c_{i+1}-N,c_i-N})$ and requiring that $A^1_{a,c_1}PA^2_{c_n-N,b} \geq   1$.

Since the set $c$ could in principle contain as many elements as the dimention of the diagram (as was the case in Figure \ref{ExampleProduct}), determining whether two nodes situated on top and bottom of the product diagram are connected, requires finding all paths of length $l$ with $2 \leq l \leq N+1$. This could be computed by raising  $A^3$ to the power $N+1$
 see Code \ref{MatrixMultiplication}

\begin{code}{\rm Below are code lines which solve the product diagram for two diagrams $D^1, D^2$ and matrix representions $A^1,A^2$ with  $dim = N$.
		\begin{equation}\label{MatrixMultiplication}
		\boxed{
		\begin{array}{l}
		 \mathrm{input} \ (A^1,A^2)\\
		 
		 \mathrm{M} = \mathrm{matrix_{3N,3N}}(x = 0)\\
		 \ \ \ \ \mathrm{M_{1:N,1:N}} =  \mathrm{A^1_{1:N,1:N}} \\
		 \ \ \ \ \mathrm{M_{1:N,(N+1):2N}} =  \mathrm{A^1_{1:N,(N+1):2N}} \\
		 \ \ \ \ \mathrm{M_{(N+1):2N,(N+1):2N}} =  \mathrm{A^1_{(N+1):2N,(N+1):2N}+A^2_{1:N,1:N}} \\
		 \ \ \ \ \mathrm{M_{(N+1):2N,(2N+1):3N}} =  \mathrm{A^2_{1:N,(N+1):2N}} \\
		\ \ \ \  \mathrm{M_{(2N+1):3N,(2N+1):3N}} =  \mathrm{A^2_{(N+1):2N,(N+1):2N}} \\
		\ \ \ \   \mathrm{M_{i,i}} =  1 \ \ \ \ \ \forall i \in [1,3N]\\
		M' = M^{N+1} \\
		(\text{finds all of the nodes which are connected to eachother by a path of length l $\leq$ N+1})\\
		
		\mathrm{ProductDiagram} = \mathrm{matrix_{2N,2N}}(x = 0)\\
		  \mathrm{for \  i \in \ 1:N}\\ 
		  \ \ \ \ x = max(which(M'_{i,:} = 1))\\
		 \ \ \ \  \mathrm{ProductDiagram}_{i,x} = 1\\
		\mathrm{for \  i \in \ (N+1):2N}\\ 
		\ \ \ \ x = max(which(M'_{i+N,:} = 1))\\
		\ \ \ \  \mathrm{ProductDiagram}_{i,x} = 1\\
		\mathrm{return \ \ \ \  ProductDiagram}
			\end{array}
		}
		\end{equation}
	} 
\end{code}
\section{Computing representations for the generators of planar diagrams} 

As mentioned, we call the $n$ dimensional vector space of planar diagrams $TL_n (d)$. For a given dimension $n$, there are $\frac{1}{n+1}$ $2n \choose n$ distinct planar diagrams which form a basis of $TL_n (d)$, with $n-1$ generators $U_1, U_2 \ldots U_{n-1}$ which are: 

\begin{equation}
 \begin{matrix}
	\begin{tikzpicture}
	\draw (0,0) rectangle (2,1);
	\draw (0.3,0)  -- (0.3,1);
	\draw (.6,0)  -- (.6,1);
	\draw (.9,0)  -- (.9,1);
     \draw (1.3,.35) node[anchor=south] {...};
   
	\draw (1.7,0)  -- (1.7,1);
	\draw (2,0)  -- (2,1);
	
	\end{tikzpicture} 
\end{matrix}\
, \ \
\begin{matrix}
\begin{tikzpicture}
\draw (0,0) rectangle (2,1);
\draw (.3,1) .. controls +(0,-0.4) and +(0,-0.4) .. (0.6,1);
\draw (.3,0) .. controls +(0,0.4) and +(0,0.4) .. (0.6,0);
\draw (.9,0)  -- (.9,1);
\draw (1.3,.35) node[anchor=south] {...};

\draw (1.7,0)  -- (1.7,1);
\draw (2,0)  -- (2,1);

\end{tikzpicture} 
\end{matrix}\
, \ \
\begin{matrix}
\begin{tikzpicture}
\draw (0,0) rectangle (2,1);
\draw (0.3,0)  -- (0.3,1);
\draw (.6,1) .. controls +(0,-0.4) and +(0,-0.4) .. (0.9,1);
\draw (.6,0) .. controls +(0,0.4) and +(0,0.4) .. (0.9,0);
\draw (1.3,.35) node[anchor=south] {...};

\draw (1.7,0)  -- (1.7,1);
\draw (2,0)  -- (2,1);

\end{tikzpicture} 
\end{matrix}
, \ \ \ldots \ \ 
\begin{matrix}
\begin{tikzpicture}
\draw (0,0) rectangle (2,1);
\draw (0.3,0)  -- (0.3,1);
\draw (0.6,0)  -- (0.6,1);
\draw (0.9,0)  -- (0.9,1);

\draw (1.3,.35) node[anchor=south] {...};
\draw (1.5,1) .. controls +(0,-0.4) and +(0,-0.4) .. (1.8,1);
\draw (1.5,0) .. controls +(0,0.4) and +(0,0.4) .. (1.8,0);

\draw (2,0)  -- (2,1);

\end{tikzpicture} 
\end{matrix}
\end{equation} 

These generators satisfy the following relations:

\begin{equation}\label{Eq:TLRelations}
\begin{array}{ll}
 & {U_i}^2=d U_i, \quad i=1,\ldots,N-1,\\
  &U_iU_{i+1}U_i=U_i,\quad i=1,\ldots,N-2,\\
  &U_iU_{i-1}U_i=U_i,  \quad i=2,\ldots,N-1,\\
  &U_iU_j=U_jU_i ,  \quad j,i=1,\ldots,N-1, \quad |j-k| \geq 2.
 \end{array}  
 \end{equation}
 
 In this section we will describe a method for finding matrix representations for the generators of $TL_n(d)$, ones which satisfy the relations $\ref{Eq:TLRelations}$. As an example, using the programming language $R$ we will find such representations for $n=4$.

\subsection{Representations of the generators for $n=4$}

\begin{figure}[h]
	\center
	\includegraphics[width=0.7\textwidth]{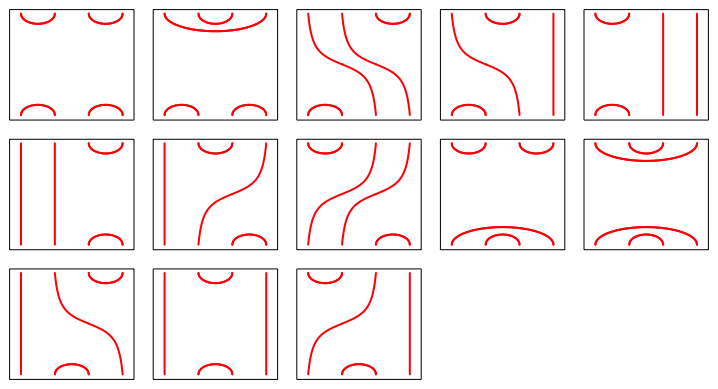}\\
	\caption{\small {All possible diagrams for the fourth dimension. We do not include the identity.} 
	}
	\label{dimension4diagrams}
\end{figure}

We begin by creating all possible diagrams in the fourth dimension, which can be seen in Figure \ref{dimension4diagrams}. We have thirteen diagrams in total and we are not including the identity. When the identity is included these diagrams form a basis for $TL_4 (d)$ and one can verify that $\frac{1}{5}$ $8 \choose 4$ is in fact $14$.

\begin{figure}[h]
	\center
	\includegraphics[height=1.5cm, width=6cm]{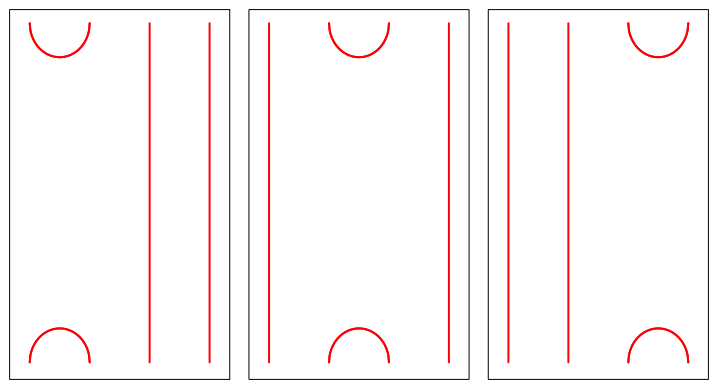}\\
	\caption{\small {The generators for the fourth dimension. Once again, we do not include the identity.} 
	}
	\label{Generators4}
\end{figure}

Next we create the generators for the fourth dimension, which can be seen in Figure \ref{Generators4}. As mentioned, these generators satisfy the relations $\ref{Eq:TLRelations}$.

$TL_4 (d)$ has ideal subspaces $I_i$ such that for any generator $U_i$ and a diagram $D_i \in I_i$ we have that $U_i D_i \in I_i$. One way of finding these ideal subspaces is by looking at the bottoms of the graphs in Figure \ref{dimension4diagrams} and grouping together all diagrams with the same bottom. We say that two diagrams have the same bottom if their bottom connections, meaning the connections from a bottom node to a bottom node, are the same. In Figure \ref{dimension4diagrams} this grouping has already been done, with the first eight diagrams belonging to the first ideal and the remaining diagrams belonging to the second ideal. Thus we can split $TL_4 (d)$ into two subspaces which, as will be shown, can be used to encode the operations of the generators by a matrix.  

\begin{figure}[h]
	\center
	\includegraphics[height=1.5cm, width=8cm]{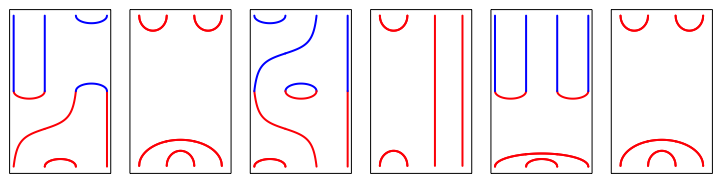}\\
	\caption{\small {Multiplication of diagrams.} 
	}
	\label{multi4}
\end{figure}   

\begin{figure}[h]
	\center
	\includegraphics[height=2cm, width=9cm]{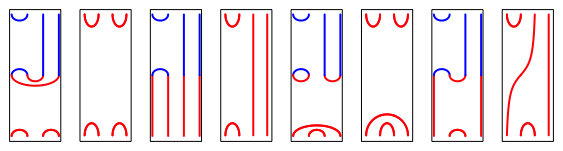}\\
	\caption{\small {Multiplication of the first generator with the second, ninth, twelfth and identity diagrams.} 
	}
	\label{multgen4}
\end{figure}

In Figure \ref{multi4}, we show three examples of the multiplication of diagrams and their outcomes. In Figure \ref{multgen4}, we multiply the first generator with some of the planar diagrams from Figure \ref{dimension4diagrams} and plot the resulting diagrams. One can verify that each diagram is mapped back into its ideal space. We use this to form a matrix representation of each generator. We form three $13 \times 13 $ matrices $\hat{U_k}$ as follows. Suppose a diagram $i$ is mapped to a diagram $j$ through multiplication by a generator $U_k$. If no $d$ is obtained by this operation then we place a $1$ in the $i$ column and $j$ row of the matrix $\hat{U_k}$. If a $d$ is obtained through this operation then we place a $d$ in the $i$ column and $j$ row of the matrix $\hat{U_k}$. We place a zero in all other components of $\hat{U_k}$. Such a matrix will represent the generator $U_k$ and will satisfy the relations $\ref{Eq:TLRelations}$.

\begin{figure}[h]
            	\center
            	\includegraphics[height=9cm, width=7.5cm]{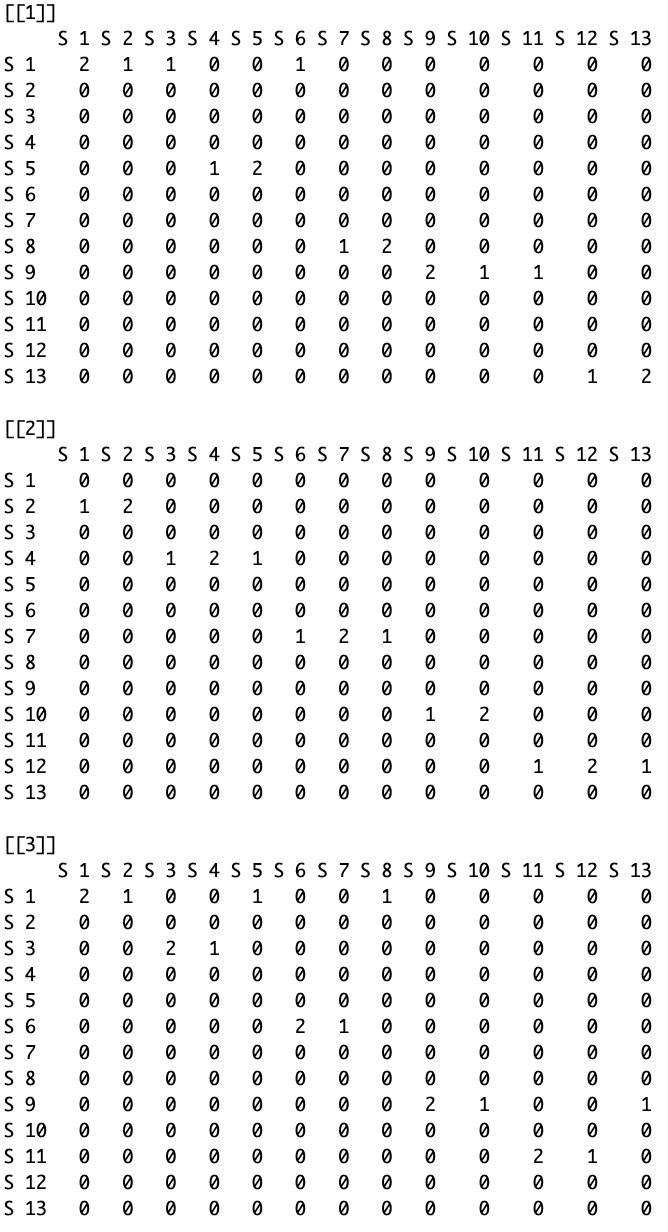}\\
            	\caption{\small {\bf } 
            	}
            	\label{representation4}
            \end{figure}

For $n=4$, we obtain the matrices given in Figure  \ref{representation4}. Instead of placing a $d$ when a bubble is created, we place a $2$ in the matrix. From Figure \ref{multgen4}, we see that the first generator times the second graph in Figure \ref{dimension4diagrams} is the first graph without a bubble. This is encoded in the first matrix, as a $1$ is placed in the second column and first row. The code also verified that the relations $\ref{Eq:TLRelations}$ were satisfied by these three matrices.

\bibliographystyle{iopart-num}

\end{document}